\documentclass[11pt]{amsart}

\usepackage[all]{xy}
\usepackage{amsmath}
\usepackage{amsfonts, mathrsfs}
\usepackage{amssymb}
\usepackage{amsthm}
\usepackage{enumerate,enumitem}
\usepackage{xcolor}
\usepackage{fullpage}
\usepackage{graphicx}
\usepackage{url}
\usepackage{hyperref}  
\usepackage{tikz} 
\usepackage{comment}
\usetikzlibrary{arrows}
\usetikzlibrary{decorations.markings}
\usepackage{changepage}
\theoremstyle{plain}
\newtheorem{theorem}{Theorem}[section]

\newtheorem{lemma}[theorem]{Lemma}

\theoremstyle{definition}

\newtheorem{example}[theorem]{Example}

\DeclareMathOperator{\inc}{inc}
\DeclareMathOperator{\wt}{Wt}

\DeclareMathOperator{\Cir}{Cir}

\title{
Quivers and moduli of their thin sincere representations in Macaulay2
}

\author{Mary Barker}
\address{Department of Mathematics and Statistics,
Washington University in St. Louis,
One Brookings Drive, St. Louis, MO, 63130}
\email{marybarker@wustl.edu}

\author{Patricio Gallardo}
\address{Department of Mathematics, 
University of California, Riverside, CA, 92521}
\email{pgallard@ucr.edu}

\date{}

\begin{document}
\begin{abstract}
We introduce the Macaulay2 package ThinSincereQuivers for studying acyclic quivers, the moduli of their thin-sincere representations, and the reflexive flow polytopes associated to them. We provide some background in the topic and illustrate how the package recovers examples from the literature.
\end{abstract}

\maketitle

\section{Introduction}
%  https://github.com/marybarker/M2
Our work, the package \emph{ThinSincere Quivers} for Macaulay2 \cite{M2}, 
provides computational tools at the crossroads of toric geometry, graph theory, and moduli spaces. 
The main objects are acyclic quivers (finite directed graphs) and representations that associate a one-dimensional vector space to each vertex, i.e., a dimension vector equal to all 1's. These representations are called 
\emph{thin sincere} ones and their moduli spaces are projective toric varieties  associated to the well-studied flow polytopes, see
 \cite{Hille1998}, \cite{Altmann1999}, \cite{hille2003quivers} and \cite{domokos2016equations}. 
Additional motivation is given by Craw and Smith \cite{Craw2008b} who showed that every toric variety is the fine moduli space for stable thin representations of an appropriate quiver with relations. 
 
Our package has many potential applications. It can be used to construct families of reflexive polytopes which play a central role in Mirror symmetry \cite{altmann2009flow}. It explicitly describes the changes among the possible flow polytopes constructed from a given quiver \cite{hille2003quivers}, and constructs  the irreducible components of the core of toric hyperk\"ahler varieties \cite{hausel2002toric}. It also provides toric compactifications of moduli spaces such as $\text{M}_{0,n}$ \cite{blume2020quivers}.  We will  illustrate such applications while recovering examples from above references. 

 The main functionality is encapsulated in the \verb|ToricQuiver| datatype. 
Quivers are equipped with either an 
integer
number for each of the arrows, i.e., an integer flow, or an integer weight for the vertices.  The type \verb|ToricQuiver| can be constructed by either an incidence  matrix and a vector of weights or a function that produces an acyclic quiver from an arbitrary graph by utilizing the Package \emph{Graphs} \cite{GraphsSource} (\S \S \ref{sec:familyQuiver}). Given a quiver, the user can study its sub-quivers and their stability
(\S\S \ref{sec:stability}), the space of regular flows and the chamber structure within the space of weights 
(\S \S \ref{Sec:GITwalls}), the flow polytope associated to a given weight (\S \S \ref{sec:flowPolytopes}), and the moduli spaces of thin-sincere representations (\S \S \ref{sec:moduli}). We also include a function for detecting if a quiver is tight (an important technical assumption) and implement an algorithm for tightening it.  

Next, we describe the structure of our package. We remark that the package documentation contains more examples which illustrate how the package can be used.

\subsection{Acknowledgement}
We would like to thank Daniel Mckenzie 
and Elo\'isa Grifo
for helpful conversations related to this work.  P. Gallardo thanks Washington University at St Louis and the University of California, Riverside for their support and welcoming environment. Mary Barker would like to thank the department of mathematics at Washington University in St. Louis for their support and access to computational resources. 
\section{Description}

\begin{figure}[h!]
\begin{tikzpicture}[label/.style args={#1#2}{%
    postaction={ decorate,
    decoration={ markings, mark=at position #1 with \node #2;}}}]
\def\ra{0.5}
\def\rc{0.5}
\begin{scope}[xshift=-4cm]
\fill (-0.2 + \rc,1) circle (1.5pt) node[left] (A1) {$v_{0}$};
\fill (4.7- \rc,1) circle (1.5pt) node[right] (A2) {$v_{1}$};
\fill (2,0) circle (1.5pt) node[right] (B1) {$v_{2}$};
\fill (2,1) circle (1.5pt) node[right] (B2) {$v_{3}$};
\fill (2,2) circle (1.5pt) node[right] (B3) {$v_{4}$};
\draw[ ultra thick, label={0.6}{[above]{$w_0$}},-> ] 
(0+ \rc,0.9) -- (1.9,0);
\draw[ ultra thick, label={0.6}{[above]{$w_1$}},-> ] (0+\rc,1) -- (1.9,1);
\draw[ ultra thick, label={0.6}{[above]{$w_2$}},-> ] (0+\rc,1.1) -- (1.9,2);
\draw[ ultra thick, label={0.6}{[above]{$w_3$}},-> ] 
(4.5-\rc,0.9) -- (2.6,0);
\draw[ ultra thick, label={0.6}{[above]{$w_4$}},-> ] 
(4.5-\rc,1) -- (2.6,1);
\draw[ ultra thick, label={0.6}{[above]{$w_5$}},-> ] 
(4.5-\rc,1.1) -- (2.6,2);
\draw (2,-1.5) node[above] (B3) 
{\texttt{BipartiteQuiver23}};
\end{scope}
%%%%%%%%%%%%%%%%%%%%%%%%%%%%%%%%%%%%%%%
%%%%%%%%%%%%%%%%%%%%%%%%%%%%%%%%%%%%%%%
%%%%%%%%%%%%%%%%%%%%%%%%%%%%%%%%%%%%%%%
\begin{scope}[xshift=2.3cm]
\fill (-0.4 +\ra,1.2) circle (1.9pt) node[left] (A1) {$v_{0}$};
\fill (4.9 -\ra -0.3 ,1.2) circle (1.9pt) node[right] (A2) {$v_{2}$};
\fill (2+0.3,0-0.3) circle (1.9pt) node[below] (B1) {$v_{1}$};
\draw[ultra thick, label={0.5}{[below]{}},-> ] 
(-0.1 + \ra,0.8) -- (1.9,-0.2);
%
% horizontal
\draw[ultra thick, label={0.5}{[above]{}}, -> ] 
(0 +\ra,1) -- (4.3 -\ra,1);
\draw[ultra thick, label={0.5}{[below]{}}, -> ] 
(0 + \ra ,1 +0.2) -- (4.3 -\ra,1+0.2);
\draw[ultra thick, label={0.5}{[below]{}}, -> ] 
(0 + \ra,1 +0.4) -- (4.3 -\ra,1+0.4);
\draw[ultra thick, label={0.7}{[above]{}},<- ] 
(4.5 -\ra,0.9-0.1) -- (2.7,0-0.1);
\draw[ultra thick, label={0.7}{[above]{}},<- ] 
(4.5+0.2 -\ra,0.9 -0.3) -- (2.6 +0.2,0-0.3);
\draw (2,-1.5) node[above] (B3) 
{\texttt{threeVertexQuiver\{1,2,3\}}};
\end{scope}
;
\begin{scope}[xshift=7cm]
\fill (2,2) circle (1.5pt) node[above] (B3) {$v_{0}$};
\fill (4,1.7) circle (1.5pt) node[above] (A2) {$v_{1}$};
\fill (4,0.3) circle (1.5pt) node[below] (A2) {$v_{2}$};
\fill (2,0) circle (1.5pt) node[below] (B1) {$v_{3}$};
\draw[ ultra thick, label={0.5}{[left]{$w_2$}},-> ] (2,1.8) -- (2,0.2);
\draw[ ultra thick, label={0.5}{[above]{$w_0$}},-> ] (2.1,2) -- (3.9,1.7);
\draw[ ultra thick, label={0.5}{[right]{$w_4$}},-> ] (4,1.5) -- (4,0.4);
\draw[ ultra thick, label={0.5}{[below]{$w_5$}},-> ] (3.8,0.3) -- (2.1,0.1);
\draw[ ultra thick, label={0.4}{[above]{$w_1$}},-> ] (2.1,1.9) -- (3.9,0.4);
\draw[ ultra thick, label={0.5}{[below]{$w_3$}},-> ] (4,1.6) -- (2.1,0.2);
\draw (3,-1.5) node[above] (B3) 
{\texttt{completeGraph 4}};
\end{scope}
\end{tikzpicture}
\caption{
Examples of acyclic quivers
%Bipartite quiver \texttt{B23} (left), three vertex quiver of type $(1,2,3)$ (center), and the quiver \texttt{K4} associated to a complete graph 4 (right), see Section \ref{sec:familyQuiver}.
}
\label{fig:example1}
\end{figure}
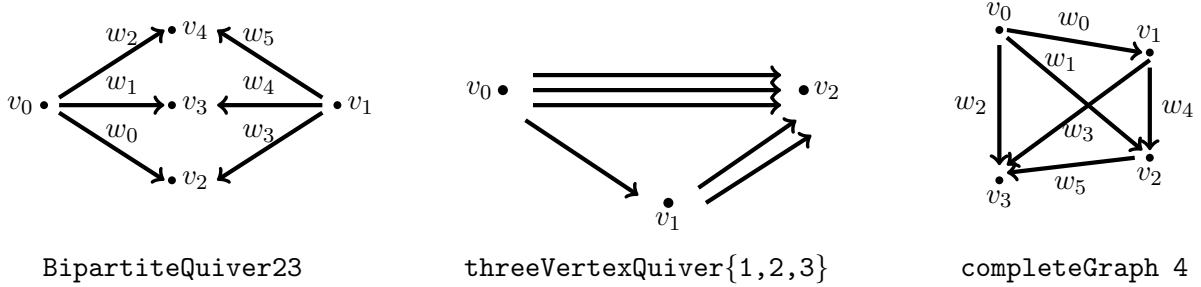

Let $Q$ be an acyclic finite directed graph i.e. an acyclic quiver. Denote the vertices by $Q_0$ and the arrows by $Q_1$, see Figure \ref{fig:example1}. We will have an integer associated to each arrow. This function $\mathbf{w}: Q_1 \to \mathbb{Z}$
is known as an \emph{integral flow} of $Q$, and it is called \emph{regular} if it is non-negative i.e., $\mathbf{w}(a) \geq 0$ for all $a \in Q_1$.
The vertices also have constants associated to them. A function $\theta: Q_0 \to \mathbb{Z}$ is called an \emph{integral weight} if $\displaystyle \sum_{i\in Q_0} \theta(i) = 0$. 
The set of all integral weights is denoted as 
$\wt(Q) \subset \mathbb{R}^{Q_0}$.
Every integral flow induces a weight on the vertices $Q_0$ as follows: Let $a\in Q_1$ be an arrow. Then,  $a^{+}\in Q_0$ will denote its head and $a^{-}\in Q_0$ will denote its tail.
The so called \emph{incidence map}  generates a weight from the flow
$\mathbf{w}$ by
\begin{equation*}
\inc(\mathbf{w})(i) := \sum_{a\in Q_1 \atop a^{+} = i} \mathbf{w}(a)  - \sum_{a\in Q_1 \atop a^{-} = i} \mathbf{w}(a) \quad \text{ for all } i \in Q_0.
\end{equation*}
In case the flow $\mathbf{w}$ is not  explicitly given, we assume that it is equal to $\mathbf{1} := (1, \ldots, 1)$.
The associated integral weight 
$\delta_Q := \inc(\mathbf{1})$ is called the \emph{canonical weight}.
If $Q$ is connected, we have the following short exact sequence:
\begin{align*}
0 \to \Cir(Q) \to \mathbb{Z}^{Q_1} \xrightarrow{\inc} \wt(Q) \to 0    \end{align*}
where the elements in $\Cir(Q)$ are called 
\emph{integral circulations}.  We remark that both the set of all integral circulations  $\Cir(Q) \subset \mathbb{Z}^{Q_1}$ and  the set of all integral weights  $\wt(Q)\subset \mathbb{Z}^{Q_0}$ forms a lattice within each space.   

\subsection{Constructing quivers}
\label{sec:familyQuiver}
To define a quiver, we can enter the edges as pairs $\{a,b\}$ where $a$ and $b$ are vertices and the arrow will take the orientation from $a$ to $b$. The user can determine a particular flow, a random flow, or use the default flow of 
$\mathbf{1} = (1,...,1)$ by leaving it undeclared.

The main datatype in this package is a type of HashTable called 
\emph{ToricQuiver}. It contains the following information:  a list $Q_0$ of integers representing the vertices, an ordered list $Q_1$ of pairs of integers, representing the edges of the quiver, and the incidence matrix representation of the underlying directed graph. This is the $|Q_0|\times|Q_1|$ matrix consisting of entries: 
\begin{align*}
a_{i,j}=
\begin{cases}
1,&\text{ if vertex $i$ is the head of the $j^{th}$ edge}
\\-1,&\text{ if vertex $i$ is the tail of the $j^{th}$ edge}
\\
0,&\text{ otherwise}
\end{cases}    
\end{align*}
In addition to these, the ToricQuiver datatype contains a list of integers representing the flow $\mathbf{w}$
and the weight $\inc(\mathbf{w})$.

\begin{example}\label{ex:defQuiver}
We generate the quiver $Q$ associated to a bipartite graph, 
see left on Figure \ref{fig:example1}, with a random flow $\mathbf{w}$ as follows: 
\begin{adjustwidth}{0.5cm}{0.5cm}{}\begin{verbatim}
i1:  loadPackage ("ThinSincereQuivers")
i2 : Q0 = {{0,2},{0,3},{0,4},{1,2},{1,3},{1,4}};
i3 : BipartiteQuiver23 = toricQuiver(Q0, Flow=>"Random")
o3 = ToricQuiver{connectivityMatrix => | -1 -1 -1 0  0  0  |           }
                                       | 0  0  0  -1 -1 -1 |
                                       | 1  0  0  1  0  0  |
                                       | 0  1  0  0  1  0  |
                                       | 0  0  1  0  0  1  |
                 flow => {71, 99, 69, 18, 46, 19}
                 Q0 => {0, 1, 2, 3, 4}
                 Q1 => {{0, 2}, {0, 3}, {0, 4}, {1, 2}, {1, 3}, {1, 4}}
\end{verbatim}\end{adjustwidth}
We change the flow of the quiver \verb|BipartiteQuiver23|
from its random value to another one \verb|w| by
\begin{adjustwidth}{0.5cm}{0.5cm}{}\begin{verbatim}
i1 : w   = {1,1,1,1,1,1};       
i2:  BPw = toricQuiver(BipartiteQuiver23, w) 
\end{verbatim}\end{adjustwidth}
\end{example}
Given a quiver $Q$ with weight $\theta \in \wt(Q)$, the user can recover its weight $\theta$, 
find a flow in the preimage 
$\text{inc}^{-1}\left( \theta \right)$, 
and detect if the underlying directed graph is acyclic by using
the commands \verb|theta(Q)|,  \verb|incInverse(th, Q)| where
\verb|th|$=\theta$, and \verb|isAcyclic(Q)| respectively. 

Additionally, given any connected graph $G=(V,E)$ generated by the Package $\emph{Graphs}$ \cite{GraphsSource}
, 
we can construct an acyclic quiver $Q$ by using the labeling of the vertices $V$. Indeed, the vertices $V$ are  indexed by integers $i \in \{0, \ldots, |V|-1 \}$. Then, given an edge joining the vertices $v_i$ and $v_j$, we set $a^{-}:=v_i$ if $i<j$. 
For example, to construct a quiver associated to the complete graph with four vertices, see  \cite[Example 2]{hille2003quivers}, we use
\begin{adjustwidth}{0.5cm}{0.5cm}{}
\begin{verbatim}
i1: needsPackage ("Graphs")
i2: CG = toricQuiver completeGraph 4
\end{verbatim}\end{adjustwidth}
Our package also has some built-in families.
The command \verb|bipartiteQuiver(r,n)| generates the quiver associated to the bipartite graph (for the case $r=2$ and $n=3$ see left in Figure \ref{fig:examples}).
The type of quiver shown at the center of Figure \ref{fig:examples} is generated by \verb|threeVertexQuiver{a,b,c}|
where \verb|a,b,c| are positive integers (see Figure \ref{fig:example1}).  
For a list \verb|L={L_1,...,L_n}| of $n$ positive integers, we define  \verb|chainQuiver {L_1,...,L_n}|  as a linear chain with $(n+1)$ vertices. The entry $L_i$ denotes the number of multiarrows between the vertices $i$ and $i+1$.  Finally, one can combine two quivers to obtain a new one by identifying either a vertex or an arrow from each. These two operations are given in the functions 
\verb|mergeOnVertex(Q1, V1, Q2, V2)| 
and
\verb|mergeOnArrow(Q1, A1, Q2, A2)|  respectively.

\subsection{Subquivers and their stability}
\label{sec:stability}

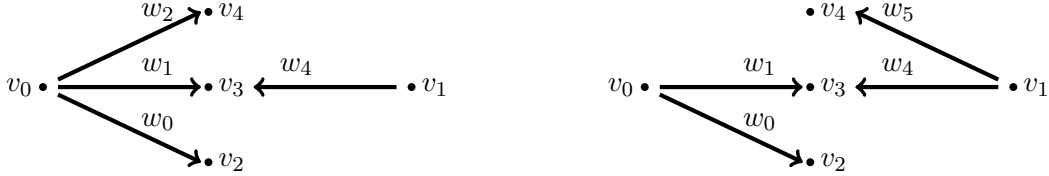
\begin{figure}
\begin{tikzpicture}[label/.style args={#1#2}{%
    postaction={ decorate,
    decoration={ markings, mark=at position #1 with \node #2;}}}]
\begin{scope}[xshift=-4cm]
\fill (-0.2,1) circle (1.5pt) node[left] (A1) {$v_{0}$};
\fill (4.7,1) circle (1.5pt) node[right] (A2) {$v_{1}$};
\fill (2,0) circle (1.5pt) node[right] (B1) {$v_{2}$};
\fill (2,1) circle (1.5pt) node[right] (B2) {$v_{3}$};
\fill (2,2) circle (1.5pt) node[right] (B3) {$v_{4}$};
\draw[ ultra thick, label={0.7}{[above]{$w_0$}},-> ] (0,0.9) -- (1.9,0);
\draw[ ultra thick, label={0.7}{[above]{$w_1$}},-> ] (0,1) -- (1.9,1);
\draw[ ultra thick, label={0.7}{[above]{$w_2$}},-> ] (0,1.1) -- (1.9,2);
%\draw[ ultra thick, label={0.7}{[above]{$w_3$}},-> ] 
%(4.5,0.9) -- (2.6,0);
\draw[ ultra thick, label={0.7}{[above]{$w_4$}},-> ] 
(4.5,1) -- (2.6,1);
%\draw[ ultra thick, label={0.7}{[above]{$w_5$}},-> ] 
%(4.5,1.1) -- (2.6,2);
\end{scope}
%%%%%%%%%%%%%%%%%%%%%%%%%%%%%%%%%%%%%%%
\begin{scope}[xshift=4cm]
\fill (-0.2,1) circle (1.5pt) node[left] (A1) {$v_{0}$};
\fill (4.7,1) circle (1.5pt) node[right] (A2) {$v_{1}$};
\fill (2,0) circle (1.5pt) node[right] (B1) {$v_{2}$};
\fill (2,1) circle (1.5pt) node[right] (B2) {$v_{3}$};
\fill (2,2) circle (1.5pt) node[right] (B3) {$v_{4}$};
\draw[ ultra thick, label={0.7}{[above]{$w_0$}},-> ] (0,0.9) -- (1.9,0);
\draw[ ultra thick, label={0.7}{[above]{$w_1$}},-> ] (0,1) -- (1.9,1);
%\draw[ ultra thick, label={0.7}{[above]{$w_2$}},-> ] (0,1.1) -- (1.9,2);
%\draw[ ultra thick, label={0.7}{[above]{$w_3$}},-> ] 
(4.5,0.9) -- (2.6,0);
\draw[ ultra thick, label={0.7}{[above]{$w_4$}},-> ] 
(4.5,1) -- (2.6,1);
\draw[ ultra thick, label={0.7}{[above]{$w_5$}},-> ] 
(4.5,1.1) -- (2.6,2);
\end{scope}
\end{tikzpicture}
\caption{
Subquivers of the bipartite quiver as defined in Example \ref{ex:defQuiver}. They are labeled by the subsets $I= \{0,1,2,4 \}$ and $\{0,1,4,5 \}$
 }
\label{fig:examples}
\end{figure}

A \emph{subquiver} $P$ of $Q$ is defined as a 
quiver $P$ such that $P_0 \subset Q_0$, $P_1 \subset Q_1$ and 
the flow function $P_1 \to \mathbb{Z}$ is obtained by restricting the one from $Q$. 
The user has two options for working with subquivers. First, they can interpret the subquiver $P \subset Q$ as one with the same vertices $P_0=Q_0$ but with a non-zero flow in the arrows $P_1 \subsetneq Q_1$. 
Since the arrows in $Q_1$ are labeled, there is a subset 
$I \subset \{0, \ldots, |Q_1|-1\}$ associated to $P_1$.
The subquiver $P$ can be constructed with the command \verb|Q^I|  
(see Figure \ref{fig:examples}). The command 
\verb|subquivers(Q)| lists all of the subquivers 
of $Q$ with the format $Q\hat{\;}I$.
The second option is to describe $P$ without any reference to $Q$.  This is done by the command \texttt{Q\_I}.  
Finally, the list of arrows defining all spanning trees can be obtained with the command
\verb|allSpanningTrees(Q)|.

Given a weight $\theta \in \wt(Q)$ and a subquiver $P\subset Q$, we can associate the concepts of $\theta$-stable,  $\theta$-semistable, and  $\theta$-unstable to $P$. They depend on a certain subset of vertices defined as follows: $V \subset Q_{0}$ is called $P$-successor closed if there is no arrow in $P_{1}$ leaving $V$.  That is, for all $a\in P_{1}$ with $a^{-}\in V$, we also have $a^{+}\in V$. 
We can check if  subset $V\subset Q_0$ is $P$-successor closed by using
\begin{adjustwidth}{.5cm}{.5cm}{}
\begin{verbatim}
i1: Q  = bipartiteQuiver(2,3)
i2: VA = {0,3}; VB = {1,3};
i3: P  = Q^{0,1,2,4} 
i4: isClosedUnderArrows(VA, P)
i5: isClosedUnderArrows(VB, P)
o4: false
o5: true
\end{verbatim}
\end{adjustwidth}
A subquiver $P \subset Q$ is $\theta$-stable (resp. $\theta$-semistable) if and only if 
$\sum_{i \in V} \theta(i) > 0$ (resp. $\geq 0$) for all successor closed $V$. 

Our package allows the user to test if a given subquiver $P$ of $Q$ is either $\theta$-stable or $\theta$-semistable by using
\begin{adjustwidth}{.5cm}{.5cm}{}
\begin{verbatim}
i1: Q = toricQuiver(bipartiteQuiver(2,3));
i2: P = Q^{0,1,4,5}; 
i3: isStable(P,Q)
i4: isSemistable(P,Q)
\end{verbatim}
\end{adjustwidth}
Notice that each quiver has always associated a weight that be can recovered with \verb|theta(Q)|. Above commands test 
stability and semi-stability with respect to it.

A quiver is {\em $\theta$-unstable} if it is not 
$\theta$-semi-stable. If $P \subset Q$ is a $\theta$-unstable quiver, then any quiver $R \subset P$ is also $\theta$-unstable. Therefore, we can define \emph{maximal $\theta$-unstable subquivers} as those that are not contained properly in any other $\theta$-unstable quiver. 
Given a quiver $Q$ with weight $\theta$, we can assume any maximal unstable quiver $P$ satisfies $P_0 = Q_0$ and $P_1 \subsetneq Q_1$.
Since $P$ is determined by a subset $I \subset Q_1$. 
We enumerate all maximal $\theta$-unstable quivers by listing such subsets. 
\begin{adjustwidth}{.5cm}{.5cm}{}
\begin{verbatim}
i1 : maximalUnstableSubquivers( bipartiteQuiver(2,3) )
o1:  HashTable{NonSingletons => {{0, 1, 2, 3}, {0, 1, 2, 4}, {0, 1, 2, 5},  
     {0, 3, 4, 5}, {1, 3, 4, 5}, {2, 3, 4, 5}}}
\end{verbatim}
\end{adjustwidth}
Similarly, if $P \subset Q$ is not $\theta$-stable, then any quiver $R \subset P$ is also not $\theta$-stable. Therefore,  we can consider 
\emph{maximal not $\theta$ stable $\theta$-subquivers}.
They are quivers that are not $\theta$-stable and are maximal with respect to the containment order. The command in this case is 
\begin{adjustwidth}{.5cm}{.5cm}{}
\begin{verbatim}
i1: maximalNonstableSubquivers( bipartiteQuiver(2,3)).
\end{verbatim}
\end{adjustwidth}
The output of above commands is in the form of a hash table. We remark that the list of unstable and not-stable quivers depends on the weight $\theta$. We illustrate such a difference for the quiver \verb|bipartiteQuiver(2,3)| but an alternative weight equal to 
$\{-5,-1,2,2,2\}$.
\begin{adjustwidth}{.5cm}{.5cm}{}
\begin{verbatim}
i1: w = incInverse({-5,-1,2,2,2},bipartiteQuiver(2,3));   -- find the flow
i2: Q = toricQuiver(bipartiteQuiver(2,3), w);             -- change the flow
i3: maximalUnstableSubquivers(Q)  
o3: HashTable{NonSingletons =>{{0, 1, 3, 4, 5},{0, 2, 3, 4, 5},{1, 2, 3, 4, 5}}}
\end{verbatim}
\end{adjustwidth}
We now turn to an important technical condition. 
A weight $\theta$ is called \emph{tight} if for every arrow  $\alpha \in Q_1$ the subquiver $P$ with $P_0= Q_0$ and 
$P_1 = Q_1 \setminus \{\alpha\}$ is $\theta$-stable. We can test this property as follows.
\begin{adjustwidth}{.5cm}{.5cm}{}\begin{verbatim}
i1: CG  = toricQuiver(completeGraph(4), {1, -2, 3, 0, 0, 0} );
i2: isTight(CG) 
o1: false
\end{verbatim}\end{adjustwidth}
The function \texttt{"makeTight"}, returns a tight quiver such that the associated flow polytope (to be defined next section) does not change. The tightening process is outlined in \cite{Altmann2009}, see also \cite{domokos2016equations}.
\begin{adjustwidth}{.5cm}{.5cm}{}\begin{verbatim}
i1: CGb = makeTight({-2,1,-2,3},toricQuiver completeGraph 4)
o1: ToricQuiver{connectivityMatrix => | -1 -1 -1 -1 | }
                                      | 1  1  1  1  |
                          flow => {-2, 1, 1, 1}
                          Q0 => {0, 1}
                          Q1 => {{0, 1}, {0, 1}, {0, 1}, {0, 1}}
                          weights => {-1, 1}

i2: isTight(CGb)
o1: true
\end{verbatim}
\end{adjustwidth}

%%%%%%%%%%%%%%%%%%%%%%%%%%%%%%%%%%%%%%%%%%%%%%%%%%%%%%%HERE 
%%%%%%%%%%%%%%%%%%%%%%%%%%%%%%%%%%%%%%%%%%%%%%%%%%%%%%%%%%%%
\section{Regular flows and cone of weights}\label{ref:SpaceRegularflows}

\begin{figure}[h!]
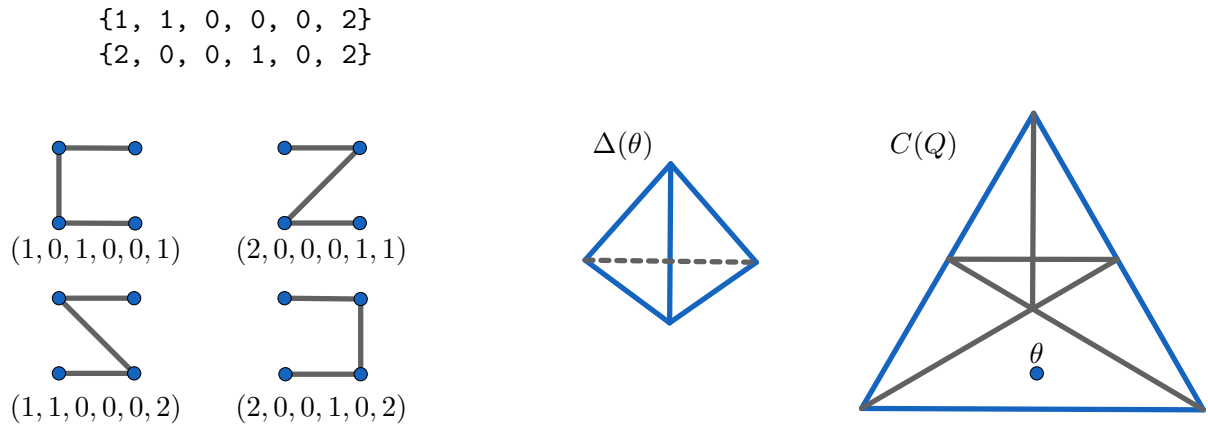

\usetikzlibrary{arrows}
\baselineskip=10pt
\hsize=6.3truein
\vsize=8.7truein
\definecolor{wrwrwr}{rgb}{0.3803921568627451,0.3803921568627451,0.3803921568627451}
\definecolor{rvwvcq}{rgb}{0.08235294117647059,0.396078431372549,0.7529411764705882}
\tikzpicture[line cap=round,line join=round,>=triangle 45,x=1.0cm,y=1.0cm]
%\fill[line width=2.pt,color=rvwvcq,fill=rvwvcq,fill opacity=0.10000000149011612] (7.714654141387821,2.0147573203060976) -- (5.458476419423932,5.959620478999879) -- (3.170213570523738,2.0332816769797573) -- cycle;
%\fill[line width=2.pt,color=rvwvcq,fill=rvwvcq,fill opacity=0.10000000149011612] (2.771611067812202,3.977666691195984) -- (1.6188382074026921,3.1726896401343687) -- (0.4949707503441151,4.006448297885045) -- (1.6332909090781587,5.287229795548277) -- cycle;
%%
%%
\begin{scope}[shift={(2.5,0)}]
\draw [line width=2.pt,color=rvwvcq] %AC 
(7.714654141387821,2.0147573203060976)-- (5.458476419423932,5.959620478999879);
\draw [line width=2.pt,color=rvwvcq] (5.458476419423932,5.959620478999879)-- (3.170213570523738,2.0332816769797573);
\draw [line width=2.pt,color=rvwvcq] (3.170213570523738,2.0332816769797573)-- (7.714654141387821,2.0147573203060976);
\draw [line width=2.pt,color=wrwrwr] (3.170213570523738,2.0332816769797573)-- (6.570246473753718,4.0157218731547);
\draw [line width=2.pt,color=wrwrwr] (4.327500060272436,4.0190233289660675)-- (7.714654141387821,2.0147573203060976);
\draw [line width=2.pt,color=wrwrwr] (5.458476419423932,5.959620478999879)-- (5.443416651842162,3.3587068392460644);
\draw [line width=2.pt,color=wrwrwr] (4.327500060272436,4.0190233289660675)-- (6.570246473753718,4.0157218731547);
\draw [fill=rvwvcq] (5.5,2.5) circle (2.5pt);
\draw (5.5,2.5) node[above]  {$\theta$};
\end{scope}  % END TRIANGLE
\begin{scope}[shift={(1.5,0)}]
\draw [line width=2.pt,color=rvwvcq] (2.771611067812202,3.977666691195984)-- (1.6188382074026921,3.1726896401343687);
\draw [line width=2.pt,color=rvwvcq]
(1.6188382074026921,3.1726896401343687)-- (0.4949707503441151,4.006448297885045);

\draw [line width=2.pt,color=rvwvcq] (0.4949707503441151,4.006448297885045)-- (1.6332909090781587,5.287229795548277);
\draw [line width=2.pt,color=rvwvcq] (1.6332909090781587,5.287229795548277)-- (2.771611067812202,3.977666691195984);
\draw [line width=2.pt,color=rvwvcq] (1.6332909090781587,5.287229795548277)-- (1.6188382074026921,3.1726896401343687);
\draw [ dashed,line width=2.pt,color=wrwrwr] (0.4949707503441151,4.006448297885045)-- (2.771611067812202,3.977666691195984);
\end{scope}
\draw [line width=2.pt,color=wrwrwr] (-2.,5.5)-- (-1.,5.5);
\draw [line width=2.pt,color=wrwrwr] (-2.,4.5)-- (-1.,5.5);
\draw [line width=2.pt,color=wrwrwr] (-2.,4.5)-- (-1.,4.5);
\draw [line width=2.pt,color=wrwrwr] (-2.,3.5)-- (-0.9882361915519591,3.4887581799061333);
\draw [line width=2.pt,color=wrwrwr] (-0.9882361915519591,2.4887581799061333)-- (-0.9882361915519591,3.4887581799061333);
\draw [line width=2.pt,color=wrwrwr] (-0.9882361915519591,2.4887581799061333)-- (-1.9882361915519597,2.4887581799061333);
\draw [line width=2.pt,color=wrwrwr] (-3.984875231619107,5.495072171252867)-- (-5.,5.5);
\draw [line width=2.pt,color=wrwrwr] (-5.,5.5)-- (-5.,4.5);
\draw [line width=2.pt,color=wrwrwr] (-3.984875231619107,4.495072171252867)-- (-5.,4.5);
\draw [line width=2.pt,color=wrwrwr] (-4.,3.5)-- (-5.,3.5);
\draw [line width=2.pt,color=wrwrwr] (-4.,2.5)-- (-5.,3.5);
\draw [line width=2.pt,color=wrwrwr] (-5.,2.5)-- (-4.,2.5);
\draw [fill=rvwvcq] (-2.,5.5) circle (2.5pt);
%\draw[color=rvwvcq] (-1.8879127819391495,5.769321707590056) node {$K$};
\draw [fill=rvwvcq] (-1.,5.5) circle (2.5pt);
%\draw[color=rvwvcq] (-0.8861910422531911,5.769321707590056) node {$L$};
\draw [fill=rvwvcq] (-2.,4.5) circle (2.5pt);
%\draw[color=rvwvcq] (-1.8879127819391495,4.761965473472907) node {$M$};
\draw [fill=rvwvcq] (-1.,4.5) circle (2.5pt);
%\draw[color=rvwvcq] (-0.8861910422531911,4.761965473472907) node {$N$};
\draw [fill=rvwvcq] (-2.,3.5) circle (2.5pt);
\draw (-1.5,1.7+2.1) node[above]  {$(2,0,0,0,1,1)$};
%
%\draw[color=rvwvcq] (-1.8879127819391495,3.7690000427002897) node {$O$};
\draw [fill=rvwvcq] (-0.9882361915519591,3.4887581799061333) circle (2.5pt);
%\draw[color=rvwvcq] (-0.8861910422531911,3.7546092393557586) node {$P$};
\draw [fill=rvwvcq] (-0.9882361915519591,2.4887581799061333) circle (2.5pt);
%\draw[color=rvwvcq] (-0.8861910422531911,2.761643808583141) node {$Q$};
\draw [fill=rvwvcq] (-1.9882361915519597,2.4887581799061333) circle (2.5pt);
%\draw[color=rvwvcq] (-1.8879127819391495,2.761643808583141) node {$R$};
\draw [fill=rvwvcq] (-3.984875231619107,5.495072171252867) circle (2.5pt);
\draw (-4.5,1.7 + 2.1) node[above]  {$(1,0,1,0,0,1)$};
%
%\draw[color=rvwvcq] (-3.8761786591946126,5.754930904245525) node {$S$};
\draw [fill=rvwvcq] (-5.,5.5) circle (2.5pt);
%\draw[color=rvwvcq] (-4.893078000997024,5.769321707590056) node {$T$};
\draw [fill=rvwvcq] (-5.,4.5) circle (2.5pt);
%\draw[color=rvwvcq] (-4.893078000997024,4.761965473472907) node {$U$};
\draw [fill=rvwvcq] (-3.984875231619107,4.495072171252867) circle (2.5pt);
%\draw[color=rvwvcq] (-3.8761786591946126,4.761965473472907) node {$V$};
\draw [fill=rvwvcq] (-4.,3.5) circle (2.5pt);
\draw (-1.5,1.7) node[above]  {$(2,0,0,1,0,2)$};
%\draw[color=rvwvcq] (-3.891356261311066,3.7690000427002897) node {$W$};
\draw [fill=rvwvcq] (-5.,3.5) circle (2.5pt);
%\draw[color=rvwvcq] (-4.893078000997024,3.7690000427002897) node {$Z$};
\draw [fill=rvwvcq] (-4.,2.5) circle (2.5pt);
%\draw[color=rvwvcq] (-3.853412256019931,2.797620816944468) node {$A_1$};
\draw [fill=rvwvcq] (-5.,2.5) circle (2.5pt);
%\draw[color=rvwvcq] (-4.85513399570589,2.797620816944468) node {$B_1$};
\draw (-4.5,1.7) node[above]  {$(1,1,0,0,0,2)$};
\draw (2.5,1.7+3.5) node[above]  {$\Delta(\theta)$};
\draw (6.5,1.7+3.5) node[above]  {$C(Q)$};
\endtikzpicture
\caption{
Cone of weights $C(Q)$ (right), polytope $\Delta(\theta)$ associated to the weight $\theta = (-2,1,-1,2)$ (center)
and $\theta$-stable trees parametrized by the vertices of the polytope
$\Delta(\theta)$ (left). Here, the quiver $Q$ is constructed from the complete graph with four vertices.
}
\label{fig:CQ}
\end{figure}

We next turn to other combinatorial structures associated to our quiver. The main reference is \cite{hille2003quivers}.
For a given weight $\theta \in \wt(Q)$, there is an integral polytope 
\begin{align*}
\Delta(\theta) =
\{  \mathbf{w} \in \mathbb{R}^{Q_1} \; | \; \inc(\mathbf{w}) = \theta  \}
\cap 
\mathbb{R}^{Q_1}_{\geq 0}
\end{align*}
which is called a \emph{flow polytope}. In general,
given a weight $\theta \in \wt(Q)$ the polytope 
$\Delta(\theta)$ can be empty. However, there exists a \emph{cone of weights} $C(Q) \subset \mathbb{R}^{Q_0}$ such that  $\Delta(\theta)$ is not empty if and only $\theta \in C(Q)$.  The polytope has the expected dimension, 
$|Q_1| - |Q_0| +1$, if and only if $\theta$ is in the interior of 
$C(Q)$. The vertices of $C(Q)$ are constructed with the so called primivite arrows from the quiver, see 
\cite[Prop. 4.7]{hille2003quivers}. The user can recover the list of primitive arrows, the cone $C(Q)$, and test the membership of $\theta$ in $C(Q)$ with the following code
\begin{adjustwidth}{0.5cm}{0.5cm}{}
\begin{verbatim}
i1: needsPackage "Polyhedra"
i2: CG =  toricQuiver completeGraph 4                         -- the quiver
i3: A = primitiveArrows CG
i4: ConeCG = coneFromVData  quiverConnectivityMatrix(CG^A)    -- C(Q)
i5: CanonicalWeight = transpose matrix {{-3,-1,2,2}};         -- weight
i5: inInterior(CanonicalWeight,ConeCG)
o1: true
\end{verbatim}
\end{adjustwidth}
We can interpret $\Delta(\theta)$ as parametrizing all possible regular flows with input equal to $\theta$. 
Each point at the interior of $\Delta(\theta)$ parametrizes a flow with entries that are strictly positive.  
The vertices of $\Delta(\theta)$ parametrize
 regular flows with input $\theta$ whose support is a spanning tree, 
 see Figure \ref{fig:CQ}.
  The user can recover all such spanning trees 
$T_i \subset Q$ by 
\begin{adjustwidth}{0.5cm}{0.5cm}{}
\begin{verbatim}
i1: CG = toricQuiver completeGraph 4;
i2: tht= {-2,1,-1,2}                   -- weight
i3: stableTrees(tht, CG)         
o1: {{4,0,5}, {0,5,3}, {0,5,2},{0,5,1}}
\end{verbatim}
\end{adjustwidth}
Each element of the output is a tuple $\{ b_1, \ldots, b_m\}$ where $b_k$ denotes the arrows of $Q$ that define the $\theta$-stable tree. In particular, we recover examples such as \cite[Fig 10]{hille2003quivers}.
Finally,  the regular flows supported on above trees can be found by using
\begin{adjustwidth}{0.5cm}{0.5cm}{}
\begin{verbatim}
i4: for i in stableTrees(tht, CG) do print incInverse(tht,CG^i)
o4: {1, 0, 1, 0, 0, 1}
    {2, 0, 0, 0, 1, 1}
    {1, 1, 0, 0, 0, 2}
    {2, 0, 0, 1, 0, 2}
\end{verbatim}
\end{adjustwidth}

%%%%%%%%%%%%%%%%%%%%%%%%%%%%%%%%%%%%%%%%%%%%%%%%%%%%%%%%%%%%%%%%%%%%

\subsection{Walls and chambers on the space of weights}\label{Sec:GITwalls}
The cone $C(Q)$ has a wall-chamber 
decomposition induced by the following equivalence relation:
$\theta \sim \theta'$ if for every subquiver $P\subset Q$, $P$ is $\theta$-semistable  implies $P$ is $\theta'$-semistable  and viceversa.
In particular, $\theta \sim t\theta$ for every $t>0$ and
$\Delta(\theta)$ is isomorphic
to $\Delta(\theta')$  up to affine lattice isomorphism. 
Our package allows the user to study this structure within $C(Q)$ and recover it for given cases.
We remark this problem with a different perspective and implementation has been solved in \cite[Sec 4]{baldoni2004counting}. For us, the starting point  to describe the walls of this decomposition is the following result due to Hille.  
\begin{lemma}
\cite[Sec 2.2]{hille2003quivers}
Each wall $W$ in $C(Q)$ is contained in a hyperplane of the form 
\begin{align*}
W(Q_0^+) = \bigg\lbrace
\theta \; \bigg| \;
\sum_{i \in Q_0^+} \theta(i) =0
\bigg\rbrace
\end{align*}
where $Q= Q_0^+ \sqcup Q_0^-$ and the full 
subquivers $Q^+$ and $Q^-$ with vertices 
$Q_0^+$ and $Q_0^-$ are connected.
\end{lemma}
We denote any wall $W(Q_o^+)$ by the subset of 
vertices  $Q_0^+$ used for defining it. The \emph{type} of the wall $W(Q_0^+)$ which is defined as  $(t^+,t^{-})$ where $t^+$ is the number of arrows starting $Q^+_0$ and ending in $Q^-_0$, and  $t^-$ is the number of arrows starting $Q^-$  and ending in $Q^+$. 
One can therefore describe a wall uniquely by the
pair:  $Q_0^+, (t^+,t^-)$. 
\begin{adjustwidth}{0.5cm}{0.5cm}{}
\begin{verbatim}
i1: L =  potentialWalls toricQuiver completeGraph 4;
i2: L#0
o1: Wall{Qplus => {1, 2, 3}}   
    WallType => (0, 3)
\end{verbatim}\end{adjustwidth}
A wall with $t^{+}$ or $t^{-}$ equal to $0$ is 
called an outer wall. Otherwise, it is called an inner wall. By \cite[Lemma 4.5]{hille2003quivers}, facets of the cone $C(Q)$ are of the form $C(Q) \cap W$ where $W$ is an outer wall.  Although every wall is contained in a hyperplane $W(Q_0)^+$, there are contiguous chambers in which the hyperplane is not a ``real wall" i.e., there exist 
two points $x$ and $y$ such that $x \sim y$ and yet $x$ and $y$ are in different half-spaces defined by the hyperplane $W$.
We can rely on the Package \emph{Lattice Polytopes} \cite{LatticePolytopesSource} for deciding whenever two polytopes are isomorphic (the required hypothesis that 
$\Delta(\theta)$ is smooth is satisfied for $\theta$ within the interior of a chamber.
\begin{adjustwidth}{0.5cm}{0.5cm}{}
\begin{verbatim}
i1: CG = toricQuiver completeGraph 4;
i1: sameChamber({-3,2,-1,2},{-2,1,-2,3}, CG)
o1: true
\end{verbatim}
\end{adjustwidth}
We remark that if the weights are not in $C(Q)$ the flow polytopes will be empty. 
\begin{adjustwidth}{0.5cm}{0.5cm}{}
\begin{verbatim}
i1: sameChamber({2,-1,1,-2},{3,-1,-1,-1}, CG)
o1: cannot be determined. stableTrees are empty
\end{verbatim}
\end{adjustwidth}
Another perspective on describing the chamber decomposition is also given by \cite{hille2003quivers}. Indeed, the spanning trees of $Q$ induce chambers within $C(Q)$ which are defined as
\begin{align*}
C_T
= \text{inc}(D_P),
&&
D_T
=
\{
\mathbf{w} \in \mathbb{R}^{Q_1}
\; | \;
\mathbf{w}(a) = 0
\text{ for all }
\alpha \notin T_1, \;
\mathbf{w}(a) \geq 0
\text{ for all }
\alpha \in T_1, 
\}.
\end{align*}
A chamber $C_T$ maybe the whole $C(Q)$, so they are not minimal with respect to the equivalence relation defined at the beginning of this section. However, the intersections  $\cap C_{T_i}$ which are of maximal dimension define a chamber system which either refines or equals the chamber decomposition of $C(Q)$, see \cite[Lemma 4.4]{hille2003quivers}. 
The relevant commands to produce the refined chamber decomposition are
\begin{adjustwidth}{.5cm}{.5cm}{}
\begin{verbatim}
i1: CG  = coneSystem(CG)       -- create the cone system
i2: rts = referenceThetas CG   -- returns a theta from each maximal chamber
o2: {{-1, -1, -1, 3}, {-2, 1, -1, 2}, {-3, 1, 1, 1}, {-2, -2, 1, 3}, 
    {-3, -1, 2, 2}, {-1, -3, 2, 2}, {-2, -2, 3, 1}}
\end{verbatim}
\end{adjustwidth}
This results recovers the chamber system for the quiver
associated to the complete graph of four vertices, see
\cite[Fig 7]{hille2003quivers}.
 A similar computation recovers the 18 chambers associated to 
 the bipartite graph with (2,3) vertices, see \cite[Sec 9]{hausel2002toric}.

\section{Applications}

We illustrate a few of the possible uses of our our packages. 
\subsection{To polytopes}\label{sec:flowPolytopes}
The polytopes $\Delta(\theta)$ are a priori contained in $\mathbb{R}^{Q_1}$ and the user can recover them formatted as a list of vertices with 
\begin{adjustwidth}{.5cm}{.5cm}{}
\begin{verbatim}
i1 : flowPolytope({-3,-3,2,2,2}, bipartiteQuiver(2,3), Format=>"Original")
o2 = {{2, 0, 1, 0, 2, 1}, {2, 1, 0, 0, 1, 2}, {0, 2, 1, 2, 0, 1}, 
     {1, 2, 0, 1, 0, 2}, {0, 1, 2, 2, 1, 0}, {1, 0, 2, 1, 2, 0}}
\end{verbatim}
\end{adjustwidth}
However, a better representation can be obtained as follows. The polytope
$\Delta(\theta)$ is contained in the fiber $\inc^{-1}(\theta)$. Therefore, we  translate them to $\Cir(Q) = \inc^{-1}(0)$.
This linear subspace is more natural because $\dim \Delta(\theta) = \dim (\Cir(Q)_{\mathbb{R}})$ for $\theta$ in the interior of $C(Q)$. 
Moreover, there is a basis of the lattice $\Cir(Q)$ for each spanning tree $T$ of $Q$. 
For example, 
given the quiver \verb|bipartiteQuiver(2,3)|
the following command computes a basis of the two dimensional vector space $\Cir_{\mathbb{R}}(Q) \subset \mathbb{R}^6$
using the spanning tree \verb|{0, 1, 4, 5}|.
\begin{adjustwidth}{.5cm}{.5cm}{}
\begin{verbatim}
i1 : basisForFlowPolytope({0, 1, 4, 5}, bipartiteQuiver(2,3))
o1 =  | 0  1  |
      | 1  -1 |
      | -1 0  |
      | 0  -1 |
      | -1 1  |
      | 1  0  |
\end{verbatim}
\end{adjustwidth}
The package automatically selects one if the none is provided and the user writes
\begin{adjustwidth}{.5cm}{.5cm}{}
\begin{verbatim}
i1: basisForFlowPolytope(bipartiteQuiver(2,3))|.
\end{verbatim}
\end{adjustwidth}
We can recover the vertices of $\Delta(\theta)$ with respect to the basis induced by any spanning tree.   For example,
if $Q$ is the bipartite graph and $\theta$ is the canonical weight. Then, $\Delta(\theta)$ is a hexagon. We obtain its vertices with respect to the basis induced by the stable tree with edges $T=\{0,1,4,5\}$.
\begin{adjustwidth}{.5cm}{.5cm}{}
\begin{verbatim}
i8 : flowPolytope({-3,-3,2,2,2}, bipartiteQuiver(2,3), Format=>{0,1,4,5})
o8 = {{0, 1}, {1, 1}, {0, -1}, {1, 0}, {-1, -1}, {-1, 0}}
\end{verbatim}\end{adjustwidth}
This output interfaces with other Macaulay2 packages such as \texttt{Polyhedra} \cite{PolyhedraSource} and
\texttt{OldPolyhedra} \cite{OldPolyhedraSource}. Then, it allows the reader to translate results based on toric Quiver computations into problems relating to polyhedral structures.  
For example, the flow polytope associated to the complete graph with 4 elements and $\theta = \{-2,1,-1,2\}$ is isomorphic to the simplex
\begin{adjustwidth}{.5cm}{.5cm}{}
\begin{verbatim}
i0: needsPackage("LatticePolytopes");
i1  Q = toricQuiver completeGraph 4
i2: w = {-2, 1, -1, 2};
i3: FromK4 = flowPolytope(w, Q);
o3: {{0, 0, 0}, {0, 0, -1}, {-1, 0, 0}, {-1, 1, 0}}
\end{verbatim}
\end{adjustwidth}
We can also explore which polytopes $\Delta(\theta)$ are reflexive. For example, given the quiver from the complete graph with four vertices and its canonical weight $\{-3,-1,1,3\}$, then we test if the associated polytope is reflexive.
\begin{adjustwidth}{.5cm}{.5cm}{}
\begin{verbatim}
i1: needsPackage ("OldPolyhedra");
i2: K4 = toricQuiver completeGraph 4;
i3: PolyK4 = convexHull transpose matrix flowPolytope({-3,-1,1,3}, K4);
i4: isReflexive PolyK4                 
o4: true
\end{verbatim}
\end{adjustwidth}
This result is expected because in fact the polytope $\Delta(\delta_Q)$ is always reflexive \cite[Prop 2.7]{Altmann1999}. 

\subsection{Applications to moduli spaces}\label{sec:moduli}
A thin sincere representation of a quiver $Q$ assigns a one dimensional vector space to each vertex in $Q_0$ and a linear map to each arrow in $Q_1$.  The space of all representations  of $Q$ is $\text{Rep}(Q) \cong \mathbb{C}^{Q_1}$.  Two representations $\omega=(w_1, \ldots, w_{|Q_1|})$ and $\omega'=(w'_1, \ldots, w'_{|Q_1|})$
are isomorphic to each other if there is an element 
$\vec v_{L} =(t_1, \ldots, t_{|Q_0|}) \in  \left( \mathbb{C}^* \right)^{Q_0}$ such that $w'_{a}t_{a^{-}} = t_{a^+}w_a$ for every arrow $a \in Q_1$.
We observe that given a representation
$ \omega$, we can define a subquiver $P := \text{Supp}(\omega)$ such that $P_0=Q_0$ and  $P_1 = \{ a \in Q_1 \;|\;  \omega(a) \neq 0 \}.$
This quiver is called the support of $\omega$. 

We say that a representation $T$ is $\theta$-stable (resp. $\theta$-semistable) if for each proper subrepresentation $S$ we have that $\sum_{q | S(q) \neq 0 }\theta(q) > 0
\; ( \geq 0 \; \text{ resp})
$. 
This notion of stability is equivalent to the one given in Section \ref{sec:stability}. Indeed, 
a representation $T$ is $\theta$-stable precisely when the subquiver 
$\text{Supp}(T)$ is $\theta$-stable, 
see  \cite[Lemma 1.4]{Hille1998} and discussion before \cite[Sec 2.2]{Altmann1999}. 
By using tools from Geometric Invariant Theory, given a quiver $Q$ with weight $\theta \in C(Q)$, there is a projective complex toric variety 
$$
\overline{M}(Q, \theta) :=
Rep(Q) / \! \! /_{\theta} (\mathbb{C}^{*})^{|Q_0|-1}
$$
of dimension $|Q_1| - |Q_0| +1$ parametrizing $\theta$-semistable thin-sincere representations up to isomorphism, see \cite{King1994} and \cite{Hille1998}.  
\begin{theorem}\cite{Hille1998}
If $\theta$ is in the interior of  $C(Q)$, the compact toric variety 
$\overline{M}(Q,\theta)$ is not empty and its associated polytope is 
equal to $\Delta(\theta)$. 
\end{theorem}
Within moduli theory it is important to describe the representations parametrized by $\overline{M}(Q,\theta)$ as well as it geometry. Above discussion implies that we can describe the $\theta$-stable and $\theta$-semistable representations parametrized by every $\overline{M}(Q,\theta)$. Moreover, via toric geometry, our package can be used to study these moduli spaces.  In particular, 
if $\theta \sim \theta'$ them  $\overline{M}(Q,\theta)$ and $\overline{M}(Q,\theta')$ are isomorphic. Therefore, the geometry of $C(Q)$ describes the birational geometry of these compactifications. Particular examples of algebraic varieties constructed via our package include the blow up of $\mathbb{P}^n$ along a linear subspace, see \cite{qin2018blow}, and toric compactifications of the moduli space of $n$ labeled points in $\mathbb{P}^1$, see \cite{blume2020quivers}.

\bibliographystyle{alpha}
\bibliography{biblio.bib}

\newcommand{\etalchar}[1]{$^{#1}$}
\begin{thebibliography}{BSDLV04}

\bibitem[AH99]{Altmann1999}
Klaus Altmann and Lutz Hille.
\newblock Strong exceptional sequences provided by quivers.
\newblock {\em Algebras and Representation Theory}, 2(1):1--17, 1999.

\bibitem[ANSW09]{altmann2009flow}
Klaus Altmann, Benjamin Nill, Sabine Schwentner, and Izolda Wiercinska.
\newblock Flow polytopes and the graph of reflexive polytopes.
\newblock {\em Discrete mathematics}, 309(16):4992--4999, 2009.

\bibitem[AVS09]{Altmann2009}
Klaus Altmann and Duco Van~Straten.
\newblock Smoothing of quiver varieties.
\newblock {\em manuscripta mathematica}, 129(2):211--230, 2009.

\bibitem[BH20]{blume2020quivers}
Mark Blume and Lutz Hille.
\newblock Quivers and moduli spaces of pointed curves of genus zero, 2020.

\bibitem[BIJ{\etalchar{+}}]{GraphsSource}
Jack Burkart, David~Cook II, Caroline Jansen, Amelia Taylor, and Augustine
  O'Keefe.
\newblock {Graphs: graphs and directed graphs (digraphs). Version~0.3.2}.
\newblock A \emph{Macaulay2} package available at
  \url{https://github.com/Macaulay2/M2/tree/master/M2/Macaulay2/packages}.

\bibitem[{B}ir]{OldPolyhedraSource}
Ren{\'e} {B}irkner.
\newblock {OldPolyhedra: {C}onvex polyhedra. Version~1.3}.
\newblock A \emph{Macaulay2} package available at
  \url{https://github.com/Macaulay2/M2/tree/master/M2/Macaulay2/packages}.

\bibitem[BLKa]{PolyhedraSource}
Ren{\'e} Birkner and {L}ars Kastner(Maintaining~author).
\newblock {Polyhedra: convex polyhedra. Version~1.10}.
\newblock A \emph{Macaulay2} package available at
  \url{https://github.com/Macaulay2/M2/tree/master/M2/Macaulay2/packages}.

\bibitem[BSDLV04]{baldoni2004counting}
Welleda Baldoni-Silva, Jes{\'u}s~A De~Loera, and Michele Vergne.
\newblock Counting integer flows in networks.
\newblock {\em Foundations of Computational Mathematics}, 4(3):277--314, 2004.

\bibitem[CS08]{Craw2008b}
Alastair Craw and Gregory~G Smith.
\newblock Projective toric varieties as fine moduli spaces of quiver
  representations.
\newblock {\em American journal of mathematics}, 130(6):1509--1534, 2008.

\bibitem[DJ16]{domokos2016equations}
M{\'a}ty{\'a}s Domokos and D{\'a}niel Jo{\'o}.
\newblock On the equations and classification of toric quiver varieties.
\newblock {\em Proceedings of the Royal Society of Edinburgh Section A:
  Mathematics}, 146(2):265--295, 2016.

\bibitem[GS]{M2}
Daniel~R. Grayson and Michael~E. Stillman.
\newblock Macaulay2, a software system for research in algebraic geometry.
\newblock Available at \url{http://www.math.uiuc.edu/Macaulay2/}.

\bibitem[Hil98]{Hille1998}
Lutz Hille.
\newblock Toric quiver varieties.
\newblock {\em Algebras and modules, II (Geiranger, 1996)}, 24:311--325, 1998.

\bibitem[Hil03]{hille2003quivers}
Lutz Hille.
\newblock Quivers, cones and polytopes.
\newblock {\em Linear algebra and its applications}, 365:215--237, 2003.

\bibitem[HS02]{hausel2002toric}
Tam{\'a}s Hausel and Bernd Sturmfels.
\newblock Toric hyperk{\"a}hler varieties.
\newblock {\em Documenta Mathematica}, 7:495--534, 2002.

\bibitem[Kin94]{King1994}
Alastair~D King.
\newblock Moduli of representations of finite dimensional algebras.
\newblock {\em The Quarterly Journal of Mathematics}, 45(4):515--530, 1994.

\bibitem[LS]{LatticePolytopesSource}
Anders Lundman and Gustav~S{\ae}d{\'e}n St{\aa}hl.
\newblock {LatticePolytopes: A \emph{Macaulay2} package. Version~1.0}.
\newblock A \emph{Macaulay2} package available at
  \url{https://github.com/Macaulay2/M2/tree/master/M2/Macaulay2/packages}.

\bibitem[Qin18]{qin2018blow}
Xuqiang Qin.
\newblock Blow ups of {$\mathbb{P}^{n}$} as quiver moduli for exceptional
  collections.
\newblock {\em arXiv preprint arXiv:1804.09544}, 2018.

\end{thebibliography}

\end{document}